\documentclass[11pt,onecolumn]{article}

 \usepackage{graphicx}

\usepackage{amssymb}
 \usepackage{amsthm}

\usepackage{amsmath}
\usepackage{amsfonts}
\usepackage{dsfont}
\usepackage[english]{babel}
\usepackage[utf8]{inputenc}
\usepackage{verbatim}

\usepackage{natbib}            
\newcommand{\be}{\begin{equation}}

\newcommand{\ee}{\end{equation}}

\newcommand{\bea}{\begin{equation*}}

\newcommand{\eea}{\end{equation*}}

\newcommand{\RU} {\mathbf{U}} 
\newcommand{\RUE}{\widehat{\mathbf U}} 
\newcommand{\ue}{\widehat{u}} 
\newcommand{\xe}{\widehat{x}} 
\newcommand{\UE}{\widehat{U}} 
\newcommand{\XE}{\widehat{X}} 
\newtheorem{theorem}{Theorem}

\newtheorem{lemma}{Lemma}

\def\N{\mathds{N}}
\def\E{\mathbb{E}}
\def\R{\mathds{R}}

\def \w{\mathrm w}
\def \q{\mathrm q}

\def \Prob{\mathrm{P}} 

\def \erfc{\text{erfc}}
\def \d{\text{MSE}}
\def \ms{\mathcal D} 
\def \mis{\mathcal S} 
\def \bP{\overline{P}}

\begin{document}




\title{Binary input reconstruction for linear systems:\\ a performance analysis}

\author{Sophie M. Fosson}
\date{}
\maketitle



\begin{abstract}
Recovering the digital input of a time-discrete linear system from its (noisy) output is a significant challenge in the fields of data transmission, deconvolution, channel equalization, and inverse modeling. A variety of algorithms have been developed for this purpose in the last decades, addressed to different models and performance/complexity requirements. In this paper, we implement a straightforward algorithm to reconstruct the binary input of a one-dimensional linear system with known probabilistic properties. Although suboptimal, this algorithm presents two main advantages: it works online (given the current output measurement, it decodes the current input bit) and has very low complexity. Moreover, we can theoretically analyze its performance: using  results on convergence of probability measures, Markov Processes, and Iterated Random Functions we evaluate its long-time behavior in terms of mean square error.

\end{abstract}
    





\section{Introduction}\label{s_intro}
Consider the input/output linear system
\be\label{system_td}\left\{\begin{array}{l} x_k=\q x_{k-1}+\w u_{k-1}~~k=1,\dots,K\\
y_k=c x_k+n_k\\
\end{array}\right.
\ee
with $K\in\N$ (possibly tending to infinity), $u_k\in\{0,1\}$ for $k=0,\dots, K-1$, $x_k\in\R$ for $k=0,\dots,K$, $y_k,n_k\in\R$ for $k=1,\dots,K$,  $\q,\w,c\in\R$, and $\q\in (0,1)$ to preserve stability. Our aim is to recover the binary input $u_k$, in an online fashion, given the output $y_k$ corrupted by a noise $n_k$. To this purpose, we retrieve a low-complexity algorithm introduced in \cite{s09} and discussed in \cite{s10,s11}, and we propose a comprehensive theoretical analysis of its performance. As a result of the analysis, we will be able to evaluate the performance  as a function of the system's parameters.

The digital signal reconstruction problem is a paradigm in data transmissions, where signals arising from finite alphabets are sent over noisy continuous channels, and in hybrid frameworks, where digital and analog signals have to be merged in the same system. 
In \cite{s09}, a particular instance of model \eqref{system_td} was derived as time discretization of a convolution system and the input estimation described  as a deconvolution problem. The same can be achieved for model \eqref{system_td}:  if we consider the system
\be\label{system}\left\{\begin{array}{l} x'(t)=ax(t)+b u(t)~~~~ t \in [0,T]\\
y(t)=cx(t)+n(t)~~~~x(0)=x_0\\
u(t), x(t), y(t),~ a, b, c \in \R, a<0
\end{array}\right.
\ee 
we have \be\label{convolution} x(t)= e^{t a}x_0+ b\int_0^t e^{a(t-s)} u(s) \mathrm{d} s.\ee
Given \be\label{bin_input}u(t)=\sum_{k=0}^{K-1} u_k \mathds{1}_{[k\tau,(k+1)\tau[}(t),~~~u_k \in \mathcal{U}=\{0,1\},~\tau>0\ee 
we can discretize in the following way: by defining
\be\begin{split}
     x_k&:=x(k \tau )~~~ \text{ for } k=0,1,\dots, K=T/\tau\\
     \q&:=e^{\tau a}\in (0,1)\\
     \w&:=b\frac{1-e^{\tau a}}{-a}=-\frac{b}{a}(1-\q)
   \end{split}\ee
we obtain
\be\label{enc}\begin{split} x_k &=\q^k x_0+b\q^k \int_0^{k\tau}e^{-as}\sum_{h=0}^{K-1} u_h
\mathds{1}_{[h\tau,(h+1)\tau[}(s)\mathrm{d} s\\
&=\q^k x_0+b\q^k\sum_{h=0}^{k-1} u_h\int_{h\tau}^{(h+1)\tau}e^{-as}\mathrm{ d}s\\
&=\q^k x_0+\frac{b}{-a}\q^k\sum_{h=0}^{k-1} u_h e^{-a(h+1)\tau}(1-e^{a\tau})\\
&=\q^k x_0+\w\sum_{h=0}^{k-1} u_h \q^{k-1-h}\\
\end{split}\ee
 from which we have the recursive formula  
\be\label{rec_encoding}   x_k=\q x_{k-1}+\w u_{k-1}. \ee 

In system \eqref{system}, recovering $u(t)$ basically consists  in the inversion of the convolution integral $ y(t)= c e^{t a}x_0+ cb\int_0^t e^{a(t-s)} u(s) \mathrm{d}+n(t)$ (where $n(t)$ represents an additive noise), which is a long-standing mathematical  ill-posed problem: small observation errors  may produce defective solutions. Several estimation approaches have been studied in the last fifty years and the literature on
deconvolution is widespread: we refer the reader to early papers \cite{tik63, tik77}  and to later  \cite{ary78,spa96,sta02,fag02},
which show also some possible applications in geophysics, astronomy, image processing and biomedical systems. For more references, see
\cite{SoThesis}.

Most of known deconvolution methods exploit the regularity of the input function to provide good estimations. This work  instead is a contribution for deconvolution in case of discontinuous input functions.

Considering a  binary alphabet, which has been chosen mainly  to keep the analysis straightforward, is  consistent with many applications: the output of several digital devices, such as computers and detection devices \cite{s10}, are binary. Nevertheless, the algorithm and the analysis presented in this paper could be generalized to larger alphabets with no much effort. 

In \cite{s09},  low-complexity decoding algorithms were introduced,  derived from the optimal BCJR \cite{BCJR} algorithm, and applied  to perform the 
deconvolution of the system \eqref{system} with $a=0$ and $b=c=1$. In this work, we apply the simplest of those algorithms, the so-called One State Algorithm (OSA for short) to the system \eqref{system_td}. We then describe the performance in terms of  Mean Square Error (MSE)
for long-time
transmissions, through a probabilistic analysis arising from the Markovian behavior of
the algorithm. The scheme of the analysis is the same proposed in \cite{s09}, but leads to completely 
different scenarios: while for  $a=0$, $b=c=1$  standard ergodic theorems for denumerable Markov Processes were sufficient to compute
the
MSE, in the present case the denumerable model does not proved the expected results, and more sophisticated
arguments are used, arising from Markov Processes, Iterated Random Functions (IRF for short) and sequences of probability measures.

The paper is organized as follows. In Section \ref{s_prob} we complete the description of the system, giving some observations and probabilistic assumptions; in Sections \ref{s_alg} and \ref{s_sim}, we present our algorithm and some simulations. The core of the paper is the performance analysis provided in Section \ref{s_analysis}. Finally we propose some concluding observations.  Notice that Sections \ref{s_prob} and \ref{s_alg} mainly retrieve the
model presented in \cite{s09} and \cite{s10}.

\subsection{Notation}
We use the following notation  throughout the paper: $\Prob$ indicates a discrete probability, while $P(\cdot,\cdot)$ is the
transition probability kernel of a Markov Process; $\E$ is the stochastic mean. 
$\mathcal{B}(\mathcal{S})$ indicates the Borel $\sigma$-field of a space $\mathcal{S}$. Given a bounded measurable function $v$ defined on a space $\mathcal{S}$, $Pv(x)=\int_{\mathcal{S}}v(y)P(x,\mathrm{d}y)$. For every measure $\mu$ on  $(\mathcal{S},\mathcal{B}(\mathcal{S}))$ and $F \in \mathcal{B}(\mathcal{S})$,  $\mu P(F)=\int_{\ms}P(x,F)\mu(\mathrm{d}x)$.
The complementary error function $\erfc$ is defined as $\erfc(x)=\frac{2}{\sqrt{\pi}}\int_x^{+\infty} e^{-t^2}\mathrm{d}t$, $x\in\R$; the indicator function $\mathds{1}_A(x)$ is equal to one if $x\in A$, and zero otherwise. Moreover we often use the following acronyms: OSA for One State Algorithm,  MSE for Mean Square Error, IRF for for Iterated Random Functions, MAP for Maximum a Posteriori.

\section{Problem Statement }\label{s_prob}
Let us develop a deeper understanding of the problem and  specify some assumptions.

First notice that, given $ x_k=\q x_{k-1}+\w u_{k-1}$, we have 
\bea x_k=\q^k x_0+\w\sum_{h=0}^{k-1}u_{k-h-1}\q^h\eea
which shows that each $x_k$ is determined by the initial state $x_0$ and by a binary  polynomial  of degree $k-1$ in $\q$.
From now onwards, let $x_0=0$, so that, for any $k=0,1,\dots,K$,
\be x_k\in\mathcal{X}=\w \left \{ \sum_{h=0}^{K} \mu_h \q^h,~\mu_h\in\{0,1\} \right\}.\ee

Moreover, let us introduce some prior probabilistic information:

\emph{Assumption 1:} the additive noise $n_k$  is white Gaussian, that is,  $n_1,\dots, n_K$  are realizations of independent Gaussian
random
variables $N_1,\dots, N_K$, with null mean and variance $\sigma^2$.

\emph{Assumption 2:} the binary inputs $u_0,\dots, u_{K-1} $ are realizations of independent Bernoulli random variables $U_0, \dots
U_{k-1}$ with parameter $\frac{1}{2}$.

Input and noise are also supposed to be mutually independent. Under these assumptions the 
system  can be rewritten in probabilistic terms as follows  (capital letters are used instead of
small letters to indicate random quantities): for $ k=1,\dots, K$,
\be
\left\{\begin{array}{l}
U_{k-1}\sim \text{Ber}\left(1/2\right)\\ 
N_k\sim \mathcal{N}(0,\sigma^2)\\ 
X_k=\q X_{k-1}+\w U_{k-1}~~(X_0=0)\\
Y_k=c X_k+N_k.\\
\end{array}\right.\ee

While Assumption 1 is realistic in physical terms, Assumption 2 is less motivated by applications, where  source bits are often not independent (for example, they may be governed by a Markov Chain). We have however imposed it for simplicity of treatment, although extensions to more sophisticated prior distributions do not require much effort. Similarly, we have chosen to propose a one-dimensional problem to make the analysis more readable, while the structure would be almost the same also for multidimensional problems (see \cite{s10}).

Given this setting, we aim at providing a method to decode the bit $u_{k-1}$ at each time step $k=1,\dots, K$, based on the current measurement  $y_k$ and  the probabilistic properties of the system. In order to perform this online recovery, the algorithm is allowed to store a few information (just one real value) about the state $x_k$.

\section{Binary Input Reconstruction }\label{s_alg}
The One State Algorithm (OSA for short) introduced in \cite{s09} fits the requirements described in the previous section.

The OSA is a suboptimal version of the Bahl, Cocke, Jelinek, and Raviv algorithm (most known as BCJR \cite{BCJR}), a prominent decoding
algorithm  used for convolutional codes. The BCJR  performs a maximum a posteriori estimation (MAP) of the input bit sequence by
evaluating the probabilities of the states of the encoder, through a forward and a backward recursion; it is optimal in the sense that
it minimizes the Mean Square Error:
\be \d(\RU,\RUE):=\frac{1}{K}\sum_{k=0}^{K-1}\mathds{E}(U_k-\UE_{k})^2 \ee where
$\RU=(U_0,\dots, U_{K-1})$ and $\RUE=(\UE_0,\dots, \UE_{K-1})$ is the estimated input sequence. In the binary case, this is equivalent
to
\bea \d(\RU,\RUE)=\frac{1}{K}\sum_{k=0}^{K-1}\Prob(U_k\neq\hat{U}_{k}).\eea
As shown in \cite{s09}, the BCJR can be adapted to the binary input deconvolution problem with optimal results, but with complexity drawbacks when
the transmission is long. This motivated the introduction of the OSA, a BCJR-based method that consists only in a  forward recursion and
that stores only one state at each iteration step. More precisely, the OSA pattern is as follows:
\begin{enumerate}
    \item Initialization of the state estimate $\xe_0$;
    \item For $k=1,\dots,K$: given $y_k \in \mathds{R}$ and $\xe_{k-1}$,
     \be\label{OSApattern}\begin{split}&\ue_{k-1}=
\left\{\begin{array}{l}
0\text{ if }  |y_k-c\q\xe_{k-1}|\leq |y_k-(c\q\xe_{k-1}+c\w)|  \\
1\text{ otherwise.}
\end{array}
\right.\\
&\xe_{k}=\q\xe_{k-1}+\w \ue_{k-1}.\end{split}\ee
\end{enumerate}
Typically, we assume to know $x_0$, so that we can initialize correctly by $\xe_0=x_0$. This point will be discussed later in Section \ref{s_analysis}. 
Given the probabilistic setting previously introduced, the OSA can also be written as:
\be\left\{\begin{array}{l}
\UE_{k-1}=\mathds{1}_{(c\q \XE_{k-1}+\frac{c\w}{2},+\infty)}(Y_k)\\
\XE_k=\q \XE_{k-1}+\w \UE_{k-1}.\\
\end{array}\right.\ee
While the BCJR estimates the probabilities of all the possible states at each step, the OSA  individuates the most likely state and
assumes it to be the correct one; on the basis of this state estimate, it decides on the current bit. The decoding is performed with
a MAP decision on the current bit, which in our probabilistic setting (Bernoulli input and Gaussian noise) reduces to the comparison between two Euclidean
distances.

The OSA is suboptimal, but presents two main good properties: (a) it is low-complexity,  both for number of computations and storage
locations; (b) it is causal, that is, it uses only the past and the present information to decode the current bit. Therefore, (a) it can
be applied to our case in which the number of states is (not countably) infinite and (b) it can be used online, making unnecessary the
complete transmission before starting deconvolution, this feature being fundamental to study long time transmissions.

In \cite{s09}, we introduced other causal algorithms: the Causal BCJR, which is a version of BCJR performing only the forward recursion,
and
the Two States Algorithm (TSA), which works basically the same as OSA, but estimates the two best states at each step with their probabilities of being the correct ones. The TSA is then oriented to soft decoding; in the differentiation
case $a=0$, $b=c=1$, it was proved to have similar performance to Causal BCJR and better than the OSA. Nevertheless,  neither the
Causal BCJR nor the TSA are efficient for system \eqref{system_td}. The first one, in fact, presents  complexity drawbacks due to the 
number of states. The second one, instead, has performance too much similar to the OSA, in spite of higher complexity: owing to the structure of the
state space $\mathcal{X}$,  the two best states turn out to be very close to each other,  which does not enhance the information
provided by the OSA.

\subsection{Similar algorithms}
We notice that our setting and decoding procedure \eqref{OSApattern} are very similar to the Decision Feedback Equalizer introduced to mitigate the effects of channels' intersymbol inference (ISI, see \cite{pulf_thesis} for a complete review). As in our case, the model considered in channel equalization is a linear system with digital input, and the goal is the input recovery for equalization purposes. Various methods have been proposed in literature and much effort has been addressed  towards complexity reduction (see, e.g., \cite{eyu88, due89, wil92, que07}). Typically, complexity is reduced by collecting information only from fixed time blocks, which is also our attempt; more precisely we consider only the current measurement and an estimation of the previous state, that is, ``minimal'' blocks of length one, but extensions to larger time blocks are possible to improve the performance.

The recovery techniques in the cited works present many analogies with ours. For example, the method introduced in \cite{que07}, if restricted to minimal blocks, differs from ours only for the introduction of a prior distribution on the state $x_k$.

Nevertheless, an outstanding difference lies in the model: channel equalization exploits the input estimate to provide a feedback equalizer to the system, while our final aim is just the input recovery.

Given the several connections, in our future work we will study possible implementations of our low-complexity algorithms, derived from BCJR, for channel equalization and propose more detailed comparisons.

\section{Simulations}\label{s_sim}

We now show  a few simulations' results, obtained by 2000 Monte Carlo Runs of 320 bit transmissions. 
\begin{figure}
   \begin{center}
   \includegraphics[width=8cm]{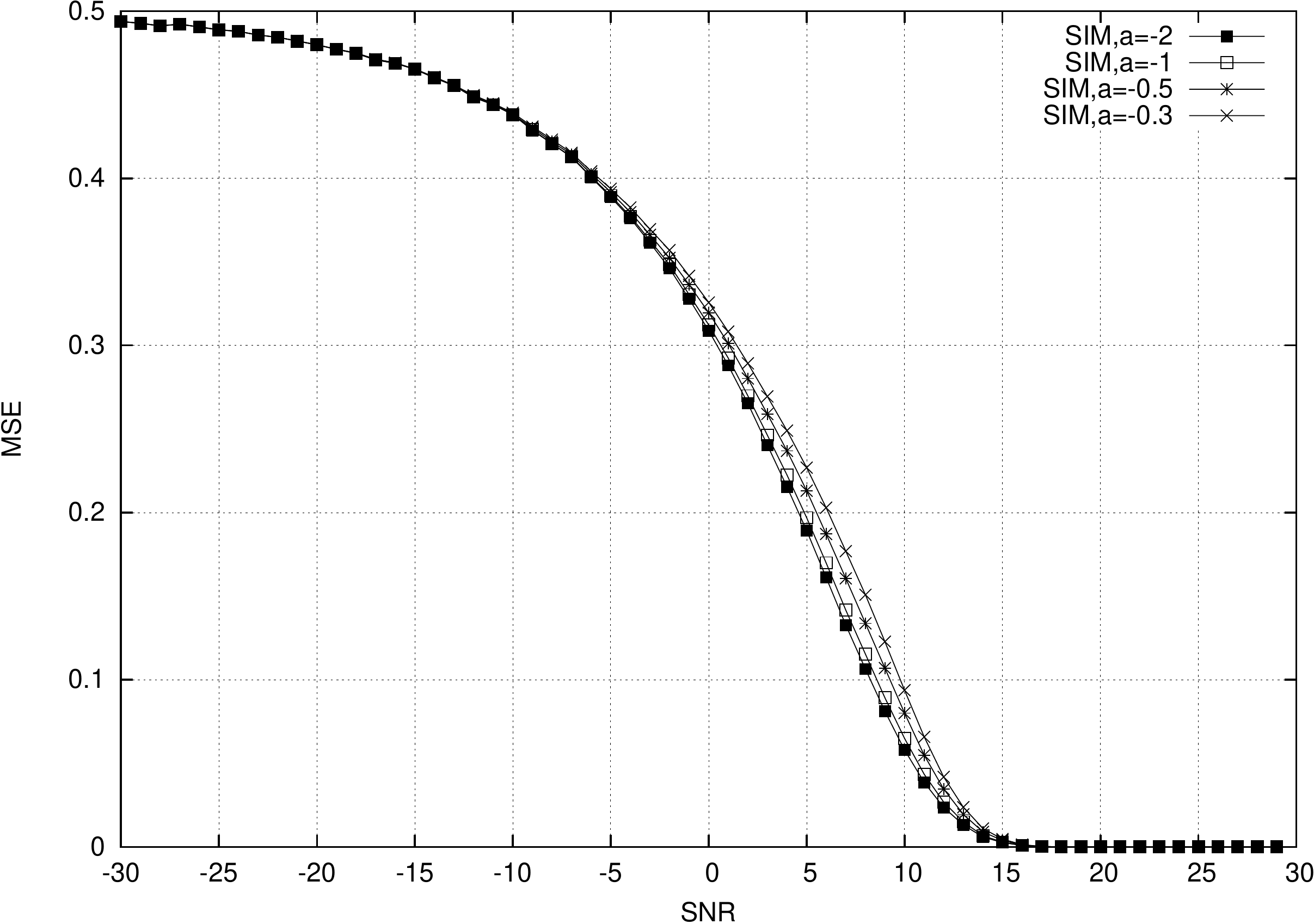}
   \caption{Mean Square Error in function of the Signal-To-Noise Ratio $\frac{c^2\w^2}{\sigma^2}$ (in dB); $b=c=1$, $a=-2,-1,-0.5,-0.3$}
\includegraphics[width=8cm]{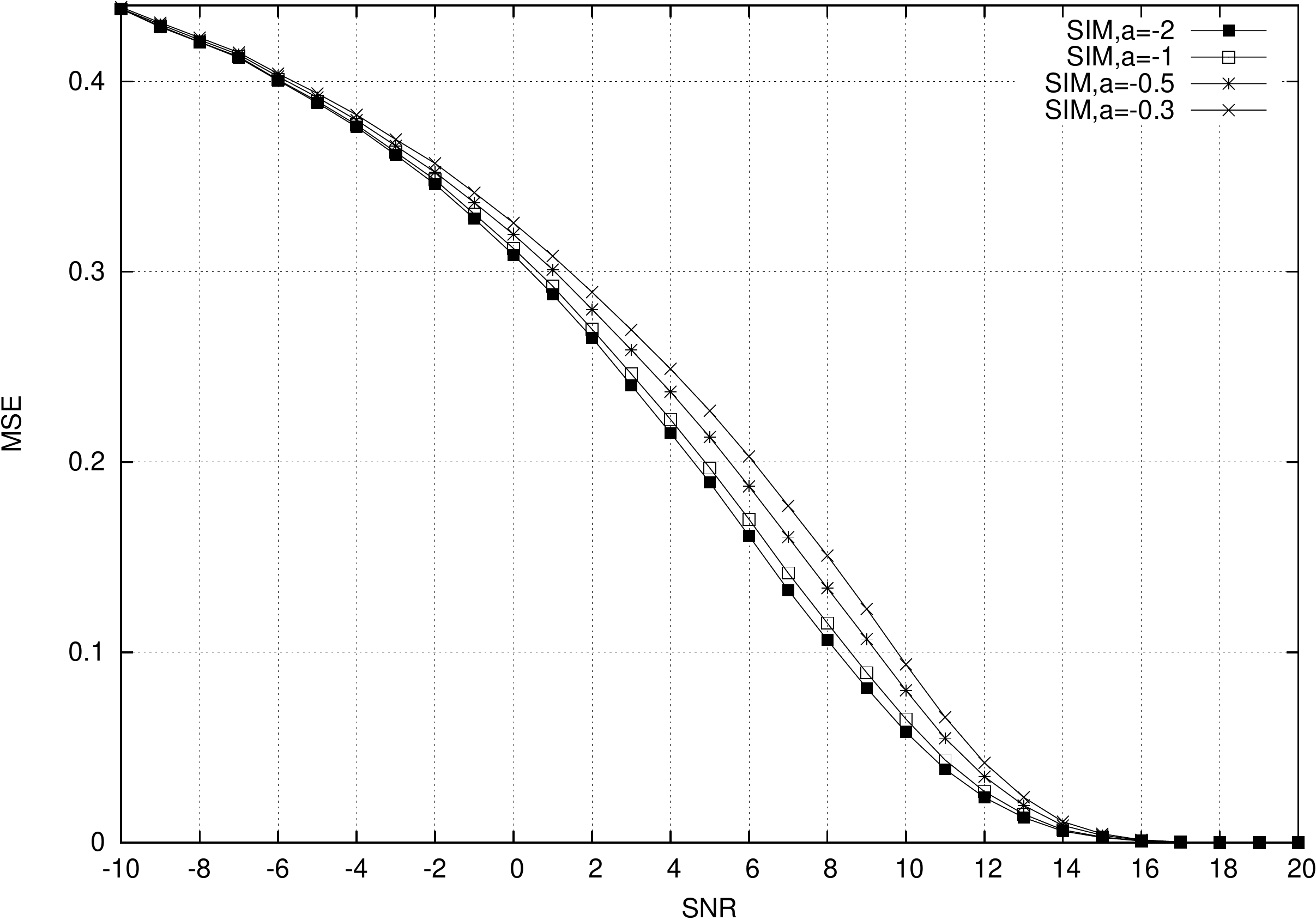}\label{figF1}
   \caption{A zoom that highlights the gain obtained by decreasing $a$.}
   \label{figF1b}\end{center}
 \end{figure} We consider the system derived from \eqref{system}, with $\tau=1$, $\q=e^{ a}$, $\w=-\frac{b}{a}(1-\q)$.  We show the behavior of the
MSE with respect to $\frac{c^2\w^2}{\sigma^2}$, that can be interpreted as the Signal-To-Noise-Ratio (SNR for short) of the
transmission: since for each $k$, the transmitted signal is $cx_k\in\{c\q
x_{k-1},c\q x_{k-1}+c\w\}$ then $c^2\w^2$ is proportional to the signal power. As expected, the MSE tends to
zero when the SNR is large, while for small SNR tends to  $\frac{1}{2}$.

If we fix $b=c=1$ and let $a$ vary, we obtain a slight gain (that is, a lower MSE curve) by decreasing $a$, as shown in Figures
\ref{figF1}-\ref{figF1b}.
In other terms, more stable systems are preferable. This phenomenon will be retrieved in Section \ref{s_analysis}.

\section{Analysis of the Algorithm}\label{s_analysis}
For simplicity, in the next we will assume $\w>0$, the analysis in the case $\w<0$ being analogous.

The goal of this section is the analytic evaluation of the Mean Square Error for the One State Algorithm, in case of long-time
transmissions.

The analysis starts from the definition of the following Markov Process:
\be\left\{\begin{array}{l}D_k=\XE_k-X_k=\q D_{k-1}+\w(\UE_{k-1}+U_{k-1})~~~~~k=1,2,\dots\\
D_0=\alpha\\
\end{array}\right.\ee
For any $k=1,2,\dots$, if  $D_{k-1}=z$ then $D_{k}\in\{\q z,\q z+\w,\q z -\w\}$, and the transition probabilities are:

\be\label{tpk}\begin{split}
P(z,\q z+\w)&=\Prob(\UE_k=1,U_k=0|D_k=z)\\
&=\frac{1}{4}\text{erfc}\left(\frac{c\q z+c\w/2}{\sigma\sqrt{2}}
\right)\\
P(z,\q z-\w)&=\Prob(\UE_k=0,U_k=1|D_k=z)\\
&=\frac{1}{4}\text{erfc}\left(\frac{-c\q z+c\w/2}{\sigma\sqrt{2}}
\right)\\
P(z,\q z)&=1-P(z,\q z+\w)-P(z,\q z-\w).\\
\end{split}\ee

Since for any $k\in\N$, $D_k\in\left \{\w  \sum_{h=0}^{k-1}
\mu_h \q^h,~\mu_h\in\{-1,0,1\} \right\}$, if we fix $D_0=0$, $(D_k)_{k\in\N}$  is a Markov Process on the denumerable state space \be\label{denumstsp}
\left \{ \w \sum_{h=0}^{\infty} \mu_h \q^h,~\mu_h\in\{-1,0,1\} \text{ such that } (\mu_h)_{h\in\N} \text{ is definitely
null}\right\}.\ee
By definitely null, we mean that for any $D_k$ the coefficients $\mu_h$ with $h\geq k$ are null. This set is denumerable since any $D_k$ can be seen as the ternary representation of  a non-negative integer. Notice that  fixing $D_0=0$ just means that $x_0$ is known.

The key point of the analysis is that, for large $k$, the MSE of the OSA can be computed using the ergodic properties of the Markov
Process $(D_k)_{k\in\N}$; more precisely, we require the existence of a stationary distribution. In the next, we will propose two different ways to study the stationary distribution: the first one does not depend on the initial state $D_0\in \ms$ (where $\ms$ is compact set which will be defined shortly), but requires some contractive properties; the second one is valid even in the non-contractive case, but depends on the initial state. 

For both methods, the presented setting is not still
adequate to study the possible stationary distributions, since the states of $(D_k)_{k\in\N}$ are transient: when the process visits a
state, then it leaves it definitely (except for a negligible set of states that have a periodic ternary representation, for example $0$, $\pm\w/(1-\q)$ ); moreover, the
process is not irreducible since there is no reciprocal communication between the states (see
\cite{s09} for more details). Thus we conclude that  no hypotheses are fulfilled to apply the
standard ergodic results for denumerable Markov
Processes  (see \cite{s09}). In other terms, if a stationary distribution exists, it does not concentrate on single states.

This suggests to consider  $(D_k)_{k\in\N}$ on a non-denumerable state space. In particular, we can
 extend  \eqref{denumstsp} to
\be\label{extstsp}\left \{ \w \sum_{h=0}^{\infty} \mu_h \q^h,~\mu_h\in\{-1,0,1\}\right\}.\ee
It is interesting to notice that if $\q\geq\frac{1}{3}$, then the set \eqref{extstsp} coincides with the closed interval
$\left[-\frac{\w}{1-\q},+\frac{\w}{1-\q}\right]$, while if $\q<\frac{1}{3}$ it is a Cantor set included in
$\left[-\frac{\w}{1-\q},+\frac{\w}{1-\q}\right]$ (for a proof of this fact, see \cite{s10}).
Let us then consider as state space
\bea \ms=\left[-\frac{\w}{1-\q},+\frac{\w}{1-\q}\right]\eea
and study the ergodic properties of the Markov Process $(D_k)_{k\in\N}$ on $\ms$. 

Before continuing the analysis, let us introduce some rigorous notions that will be used in the next. 

Let $\mathcal{B}(\ms)$ be the Borel $\sigma$-field of $\ms$.  We call  \textit{transition probability kernel} (see, e.g.,
\cite[Section 3.4.1]{mey93})  an application $P:\ms \times \mathcal{B}(\ms)\to[0,1]$ such that\\ 
(i) for each $F \in\mathcal{B}(\ms)$, $P(\cdot, F)$ is a non-negative measurable function;\\
(ii) for each $x \in \ms$, $P(x, \cdot)$ is a probability measure  on $(\ms,\mathcal{B}(\ms))$.

Given a bounded measurable function $v$ on $\ms$, we denote by $Pv$ the bounded measurable function 
defined as
\bea
Pv(x)=\int_{\ms}v(y)P(x,\mathrm{d}y).\eea
Furthermore, let $\mu$ be a measure on  $(\ms,\mathcal{B}(\ms))$: we define the measure $\mu P$ by 
\bea\mu P(F)=\int_{\ms}P(x,F)\mu(\mathrm{d}x)\;\;\;\;\;F \in \mathcal{B}(\ms).\eea
We finally define the $n$-th power of the transition kernel $P$ by $P^1(x,F)=P(x,F)$ and 
$P^n(x,F)=\int_{\ms}P(x,\mathrm{d}y)P^{n-1}(y,F)$. It is easy to see that $P^n(x,F)$ are transition kernels, too.

At this point, we can make explicit the relationship between the MSE and $(D_k)_{k\in\N}$. In the next, we will always consider $D_0=0$, if not differently specified.  Given the transition probability kernel $P$
of $(D_k)_{k\in\N}$, defined by \eqref{tpk}, we have
\be\begin{split}\label{acheserveilMP}\d&(\RU,\RUE) =\frac{1}{K}\sum_{k=0}^{K-1}\Prob(\UE_k\neq U_k)=\\
&=\frac{1}{K}\sum_{k=0}^{K-1}\int_{\mathcal{D}}\Prob(\UE_k\neq
U_k|D_k=z)P^k(\alpha,\mathrm d z)=\frac{1}{K}\sum_{k=0}^{K-1}P^k g(\alpha)
\end{split}\ee
where \be\label{def_g}g(z)=\Prob(\UE_k\neq U_k|D_k=z)\ee
and $D_0=\alpha$ is any initial state in $\ms$. Therefore $P^k g$ (and in particular its behavior for large $k$) will be the object of our further
analysis.

In the sequel, we will distinguish two main scenarios: when $(D_k)_{k\in\N}$ has some
\emph{contractive} properties and when it has not. In the first scenario, we can exploit the theory of Iterated Random Functions to
prove that $P^k g$  converges, while in the second one we will use known results of convergence of probability measures.

\subsection{Contractive case}
Let $ l\text{-Lip}(\ms)$ be  the set of all the Lipschitz functions with Lipschitz constant equal to $l$ on $\ms$. We define the
Kantorovich (or Wasserstein) distance $d_W$  
between probability measures (see \cite[Section 2.1, Example 3.2.2]{rac91}) as
\be\label{Kantorovich} d_W(\mu,\nu)=\sup_{f\in 1\text{-Lip}(\ms)}\bigg|\int_{\ms} f\mathrm d (\mu- \nu) \bigg|.\ee

We can prove the following
\begin{theorem}\label{ERGT}
If $\frac{c^2\w^2}{\sigma^2}>4$ and $\q <\frac{1}{3+\sqrt{\frac{2}{e\pi}}}$  or if $\frac{c^2\w^2}{\sigma^2}\leq 4$ and
$\q\leq \frac{1}{1+\sqrt{\frac{2}{e\pi}}}$,
then,
\be\label{ERG} \lim_{K\to\infty}\d(\RU,\RUE)=\int_{\mathcal{D}}g \mathrm  d\mu~~~\text{ for any } D_0\in \ms\ee
where $\mu$
is the unique probability measure such that  $\sup_{x\in\ms} d_W(P^k(x,\cdot),\mu(\cdot)) \stackrel{k\to\infty}{\longrightarrow}0$,
 $P$ being the kernel of $(D_k)_{k\in \N}$.\end{theorem}

Notice that $g(z)$ is time-invariant (i.e., does not depend on $k$) and can be analytically computed. In fact, given $D_{k-1}=z$,
$D_k=\q z$ if and only if $\UE_k=U_k$,  $D_k=\q z+\w$ if and only if $\UE_k=1$ and $U_k=0$, $D_k=\q z-\w$ if and only if $\UE_k=0$ and
$U_k=1$ and
\be g(z)=P(z,\q z+\w)+P(z,\q z -\w).\ee
Furthermore, the probability measure $\mu$ can be numerically evaluated.

Recall  that, as already noticed, $\frac{c^2\w^2}{\sigma^2}$ can be interpreted as the Signal-To-Noise-Ratio. Moreover, the bounds on $\q$ can be interpreted as the necessity of a stronger stability for convergence.


In order to prove the theorem, let us introduce some elements from the Iterated Random Functions theory.
\subsubsection{Iterated Random Functions}\label{IRF}
Let $(\ms, d)$ be a complete metric space and $\mis$ be a measurable space. Consider a measurable function $w: \ms \times
\mis \to \ms$ and for each fixed $s \in \mis$, $w_s(x):=w(x,s)$, $x\in\ms$. Let $(I_k)_{k\in\N}$ be a stochastic sequence in $\mis$ such
that $I_0, I_1, \dots$ are independent, identically distributed. Then, the set $\{w_{I_k}(x),~k\in\N\}$ is a family of random
functions. 
The systems obtained by iterating such random functions, called Iterated Random Functions (IRF), are studied for diverse
purposes: for example, IRF with contractive properties are used to construct fractal sets, see \cite{hut81, dia99}. More interesting for
our study is the exploitation of IRF to study Markov Processes. Given an IRF and a starting state $x\in\ms$, we can  define the induced
Markov Process  $(Z_k(x))_{k\in\N}$ as 
\be Z_k(x):=w_{I_{k-1}}\circ w_{I_{k-2}}\circ w_{I_{k-3}}\circ\cdots\circ w_{I_{0}}(x)~~~~(k\geq 1)\ee
and analyze its asymptotic behavior through the properties of $w_{I_k}(x), k\in\N$. It has been proved that if the $w_{I_k}(x)$ have
\emph{some} contractive properties, the transition probability kernel of $Z_n(x)$  converges to a probability measure, unique for all
the initial states
$x\in \ms$. The required contractive properties may be slightly different: \cite{dia99} studied the case of Lipschitz functions
$w_{I_k}(x)$ ``contracting on average'', while similar results have been obtained by \cite{sten01} without the continuity requirement on
$w_{I_k}(x)$,  by \cite{stein99} for ``locally contractive'' functions, and  by \cite{jar01} for ``non-separating on average''
functions. A useful survey on the argument has been recently proposed by \cite{ios09}.

Let us show how to exploit the IRF theory in our framework.

The evolution of $(D_k)_{k\in\N}$ can be modeled by IRF. We consider the complete metric space $\ms=\left[-\frac{\w}{1-\q},+\frac{\w}{1-\q}\right]$ naturally endowed with the
Euclidean metric $d$ of $\R$, the measurable space $\mis=\{0,1\}\times \R$ and the stochastic process $$I_k=(U_k,N_{k+1}),~~k\in\N$$ on
$\mis$, and we define the random function
\be w_{I_k}(x)=\q x+\w \mathds{1}_{(c\q x+c\w\left(\frac{1}{2}- U_k\right),+\infty)}(N_{k+1})-\w U_k,~~x\in\ms\ee that describes
the
dynamics of $(D_k)_{k\in\N}$. The key result for our purpose is the following theorem (here stated for compact spaces), which does not
require continuity:

\textbf{Stenflo's Theorem} \cite[Theorem 1]{sten01}
Suppose that there exists a constant $l<1$ such that \be\label{av_conv} \E[d(w_{I_0}(x),w_{I_0}(y))]\leq l~d(x,y)\ee for all
$x,y\in\ms$, $(\ms,d)$ being a compact metric space.
Then there exist a unique probability measure $\mu$ and a positive constant $\gamma_{\ms}$ such that
\be \sup_{x\in\ms} d_W(P^n(x,\cdot),\mu(\cdot))\leq \frac{\gamma_{\ms}}{1-l}l^n~~~~n\geq 0\ee
where $P^n(x,\cdot)$ is the $n$-step transition probability kernel of the Markov Process $Z_n(x)$.\\

Now, Theorem \ref{ERGT} can be proved by applying the Stenflo's Theorem.

\subsubsection{ Proof of Theorem \ref{ERGT}}
Let us analyze the condition (\ref{av_conv}). Consider $x,y\in\ms$ with $x>y$  (recall that $\q>0, \w>0$). Let $ H=H(x,y,I_0)$ and
$\mathcal{I}_u$ respectively be defined as
\bea\begin{split}& H:= \mathds{1}_{\left(c\q y+c\w\left(\frac{1}{2}- U_0\right),c\q x+c\w\left(\frac{1}{2}-
U_0\right)\right)}(N_1)\\ 
&\mathcal{I}_u:=\frac{1}{\sqrt{2\pi}\sigma}\int_{c\q y+c\w\left(\frac{1}{2}-
u\right)}^{c\q x+c\w\left(\frac{1}{2}- u\right)}e^{-\frac{n^2}{2\sigma^2}}\mathrm d n\\  
&=\frac{1}{2}\erfc\left(c\frac{\q y+\w\left(\frac{1}{2}-
u\right)}{\sigma\sqrt{2}}\right)-\frac{1}{2}\erfc\left(c\frac{\q x+\w \left(\frac{1}{2}-
u\right)}{\sigma\sqrt{2}}\right).\\
\end{split}\eea
Hence,
\be\label{my_contr}
\begin{split}
&\E\left[|w_{(U_0,N_1)}(x)-w_{(U_0,N_1)}(y)|\right]
=\E\left[ |\q(x-y)-\w H | \right]\\
&=\sum_{u\in\{0,1\}}\Prob(U_0=u) \frac{1}{\sqrt{2\pi}\sigma}\int_{\R}e^{-\frac{n^2}{2\sigma^2}}|\q(x-y)-\w H | \mathrm{d} n\\
&=\frac{1}{2}\sum_{u\in\{0,1\}}\frac{1}{\sqrt{2\pi}\sigma}\int_{\R}e^{-\frac{n^2}{2\sigma^2}}|\q(x-y)-\w H |\mathrm{d} n\\
&=\frac{1}{2}\sum_{u\in\{0,1\}}|\q(x-y)-\w|\mathcal{I}_u+\q(x-y)(1-\mathcal{I}_u).\\
\end{split}\ee 
If $\q(x-y)>\w$, then $\E\left[|w_{(U_0,N_1)}(x)-w_{(U_0,N_1)}(y)|\right]<\q(x-y)$ and the contraction would be proved with $l=\q$.
This is never the case when $\q<\frac{1}{3}$ and  $|x-y|\leq 2\frac{\w}{1-\q}< \frac{\w}{\q}$.

Let us then consider  $\q(x-y)<\w$. We can write
\be\begin{split}\E&\left[|w_{(U_0,N_1)}(x)-w_{(U_0,N_1)}(y)|\right]\\
&=\frac{1}{2}\sum_{u\in\{0,1\}}(\w-\q(x-y))\mathcal{I}_u+\q(x-y)(1-\mathcal{I}_u)\\
&\leq\frac{1}{2}\sum_{u\in\{0,1\}}\w\mathcal{I}_u+\q(x-y).\\
\end{split}\ee
The last expression is obtained by neglecting $-\sum_{u\in\{0,1\}}~ \q$ $(x-y)\mathcal{I}_u$, which is the sum of two second degree
terms in $(x-y)$, since, by the integral mean value theorem,
 \be \mathcal{I}_u= \frac{1}{\sqrt{2\pi}\sigma}c\q (x-y)e^{-\frac{n_0^2}{2\sigma^2}}\ee for some $n_0\in\left(c\q
y+c\w\left(\frac{1}{2}-u\right),c\q x+c\w\left(\frac{1}{2}-u\right) \right)$, $(n_0 \neq 0).$ 
The remaining terms are of order one.
 Notice  that
\be\frac{1}{2}\sum_{u\in\{0,1\}}\w\mathcal{I}_u+\q(x-y)=F(x)-F(y) \ee
where \bea F(x)=\q x-\frac{\w}{4}\erfc\left(\frac{c\q x+c\frac{\w}{2}}{\sigma\sqrt{2}}\right)-\frac{\w}{4}\erfc\left(\frac{c\q
x-c\frac{\w}{2}}{\sigma\sqrt{2}}\right).\eea
Therefore, the thesis is achieved if $F(x)$ is a contraction; since $F(x)$ is differentiable and monotone increasing, its Lipschitz
constant is the maximum of its first derivative:\be\label{F1}\begin{split}&
F'(x)=\q+\\&\frac{c\w\q}{2\sigma\sqrt{2\pi}}\left[\exp\left(-\frac{\left(c\q x+c\frac{\w}{2}\right)^2}{2\sigma^2}\right)
+\exp\left(-\frac{\left(c\q x-c\frac{\w}{2}\right)^2}{2\sigma^2}\right)\right].\\\end{split}\ee
In order to find the maximum of $F'(x)$, let us compute $F''(x)$:
\bea\begin{split}&F''(x)=\\ &=-\frac{c\w\q}{2\sigma\sqrt{2\pi}}\frac{2\left(c\q x+c\frac{\w}{2}\right)c\q}{2\sigma^2}
\exp\left(-\frac{\left(c\q x+c\frac{\w}{2}\right)^2}{2\sigma^2}\right)+\\
&~~~-\frac{c\w\q}{2\sigma\sqrt{2\pi}}\frac{2\left(c\q x-c\frac{\w}{2}\right)c\q}{2\sigma^2} \exp\left(-\frac{\left(c\q
x-c\frac{\w}{2}\right)^2}{2\sigma^2}\right)\\
\end{split}\eea
which is null for  $x$ satisfying:
\be\label{exp_eq} \left(\q x+\frac{\w}{2}\right)\exp\left(-\frac{c^2\q\w x}{\sigma^2}\right)+ \left(\q x-\frac{\w}{2}\right)=0 \ee
a solution of which is $x=0$. Now, considering that $F'(x)$ is determined by a mixture of two Gaussians, two cases may occur: (a) $x=0$
is the maximum of $F'(x)$; (b) $x=0$ is a minimum for $F'(x)$ and there are two symmetric maxima ($F''(x)$ is an even function) at
$x_0\in (0,\frac{\w}{1-\q}]$ and $-x_0$, but $x_0$ cannot be analytically computed from the exponential equation (\ref{exp_eq}). Let us study the sign of $F''(x)$ for $x\to 0$ in order to determine the nature of the point $x=0$ for $F'(x)$. Notice that 
\be\begin{split}&F''(x)>0 \Leftrightarrow 
-\left(\q x+\frac{\w}{2}\right)\exp\left(-\frac{c^2\q\w x}{\sigma^2}\right)- \left(\q x-\frac{\w}{2}\right)>0.\\
\end{split}\ee
Moreover, if $x\to 0$, 
$\exp\left(-\frac{c^2\q\w x}{\sigma^2}\right)\sim 1 - \frac{c^2\q\w x}{\sigma^2} $ and 
\be\begin{split}&-\left(\q x+\frac{\w}{2}\right)\exp\left(-\frac{c^2\q\w x}{\sigma^2}\right)- \left(\q x-\frac{\w}{2}\right)\sim-2\q x+\frac{c^2\q\w^2 x}{2\sigma^2}.\\
\end{split}\ee
Finally, if $ \frac{c^2\w^2}{\sigma^2}>4$
\be -2\q x+\frac{c^2\q\w^2 x}{2\sigma^2}\to 0^+ \text{ for } x\to 0^+\ee
\be -2\q x+\frac{c^2\q\w^2 x}{2\sigma^2} \to 0^- \text{ for } x\to 0^-.\ee
In conclusion, $x=0$ is a maximum point for $F'(x)$ if and only if  $\frac{c^2\w^2}{\sigma^2}<4$, that is, only for large noise.

Let us now study $\frac{c^2\w^2}{\sigma^2}>4$ ($x=0$ is a minimum point) and let us state conditions that make $F(x)$  contractive.
In particular, consider
$x>0$ and $\sigma^2$ close to zero: by (\ref{exp_eq}),  $|x-\frac{\w}{2\q}|$ tends to zero more quickly than $\sigma^2$, hence
$\exp\left(-\frac{\left(c\q x-c\frac{\w}{2}\right)^2}{2\sigma^2}\right)$ tends to one and the maximum of $F'(x)$ (see \eqref{F1}) may
assume very large values.

More in general, we observe that the points $x=\pm\frac{\w}{2\q}$ are tricky as they are the unique points where the OSA fails: for
these values, the error probability given by (\ref{tpk}) is $\frac{1}{2}$, no matter which is the noise variance. This ``singular''
phenomenon is more evident when the noise is small; in terms of $F(x)$, it causes large variations (then the loss of the contractivity)
in a neighborhood of the point  $\pm\frac{\w}{2\q}$, the radius of the neighborhood being larger for smaller $\sigma^2$.

Let us set in the case  $\q<\frac{1}{3}$, so that $\pm\frac{\w}{2\q}\notin \ms$. Under this
assumption, for any $x\in\ms$,
\be\begin{split}
&\exp\left(-\frac{\left(c\q
x+c\frac{\w}{2}\right)^2}{2\sigma^2}\right)<\exp\left(-\frac{\left(-c\q\frac{\w}{1-\q}+c\frac{\w}{2}\right)^2}{2\sigma^2}\right)\\
&\exp\left(-\frac{\left(c\q x-c\frac{\w}{2}\right)^2}{2\sigma^2}\right)<\exp\left(-\frac{\left(c\q
\frac{\w}{1-\q}-c\frac{\w}{2}\right)^2}{2\sigma^2}\right)\\
\end{split}\ee
hence
\be F'(x)\leq\q+\frac{c\w\q}{\sigma\sqrt{2\pi}}
\exp\left(-\frac{ c^2 \w^2\left(\frac{1-3\q}{2(1-\q)}\right)^2}{2\sigma^2}\right).\ee
As for $t\geq 0$, $\max te^{-t^2}=\frac{1}{\sqrt{2e}}$,
\bea \frac{\q}{\sqrt{\pi}}\frac{2}{1-3\q}\frac{1}{\sqrt{2e}}<1 \Longrightarrow F'(x)<1\eea

In conclusion a sufficient condition for average contractivity is

\bea\q<\frac{1}{3+\sqrt{\frac{2}{e\pi}}}.\eea

Let us now study the case $\frac{c^2\w^2}{\sigma^2}<4$ ($x=0$ is the maximum point):
\begin{equation}
\begin{split}
F'(x)&\leq F'(0)=\q+2\frac{c\w\q}{2\sigma \sqrt{2\pi}}\exp\left(-\frac{c^2\w^2}{8\sigma^2}\right) \\
       &\leq\q+2\frac{\q}{ \sqrt{2e\pi}}\\  
\end{split}
\end{equation}
then $F'(x)<1$ when $$\q<\frac{1}{1+\sqrt{\frac{2}{e\pi}}}.$$

In conclusion, we have stated that  if $\frac{c^2\w^2}{\sigma^2}>4$ and $\q<\frac{1}{3+\sqrt{\frac{2}{e\pi}}}$  or if 
$\frac{c^2\w^2}{\sigma^2}\leq 4$ and
$\q\leq \frac{1}{1+\sqrt{\frac{2}{e\pi}}}$, then the hypotheses of Stenflo's Theorem are fulfilled. 

Now, let us prove the
convergence of the Mean Square Error. 

Since $g\in L_g$-Lip$(\ms)$ where $L_g=\max_{z\in\ms}
|g'(z)|$ and 
\be |g'(z)|=\frac{c\q}{2\sigma\sqrt{2\pi}}\left|-e^{-\frac{(c\q z+c\w)^2}{2\sigma^2}}
+e^{-\frac{(-c\q z+c\w)^2}{2\sigma^2}}\right|\leq\frac{c\q}{\sigma\sqrt{2\pi}} \ee
then $L_g\leq \frac{c\q}{\sigma\sqrt{2\pi}}$ is finite. Since for any $L>0$ 
 \bea
 \begin{split}
&\sup_{f\in L\text{-Lip}(\ms)}\bigg|\int_{\ms} f \mathrm d (\mu- \nu)\bigg| = \sup_{f\in L\text{-Lip}(\ms)}L\bigg|\int_{\ms}\frac{1}{L}
f \mathrm d
(\mu- \nu)\bigg|\\ &\leq  \sup_{f\in 1\text{-Lip}(\ms)}L\bigg|\int_{\ms} f \mathrm d (\mu- \nu)\bigg|= L~ d_W(\mu,\nu)   
 \end{split}
\eea
we have 
\be
\begin{split}
  \sup_{x\in\ms}\bigg| P^k g(x)&-\int_{\ms} g\mathrm d \mu\bigg| =\sup_{x\in\ms}\bigg | \int_{\ms} g(z)P^k(x,\mathrm d z)-\int
_{\ms} g\mathrm d
\mu\bigg|
\\ &\leq \sup_{x\in\ms} \sup_{f\in L_g\text{-Lip}}\bigg|\int_{\ms} f(z)P^k(x,\mathrm d z)-\int_{\ms} f\mathrm d \mu\bigg|\\  
&= \sup_{x\in\ms} L_g d_W(P^k(x,\cdot), \mu(\cdot))\stackrel{k\to\infty}{\longrightarrow}0.\\  
\end{split}
\ee
The convergence is then assured also for the Cesàro sum, for any initial state $D_0=\alpha\in\ms$:
\be \frac{1}{K}\sum_{k=0}^{K-1}P^k g(\alpha)\stackrel{K\to\infty}{\longrightarrow}\int_{\ms} g\mathrm d \mu~~~~\forall \alpha\in\ms.\ee
\qed

Notice that the initial state does not affect the convergence value if it is contained in $\ms$, but $\ms$ has been obtained by fixing $D_0=0$: this seems not coherent. However, even if we consider $D_0=\alpha\notin \ms$, given the dynamics of $D_k$ ($\alpha$ is multiplied by $\q$ at each step), $\ms$ turns out to be the ``limit'' state space, and with high probability $D_k$ enters $\ms$ for some finite $k$, so it makes sense to reduce to $\ms$. Further details are here omitted for brevity, but one must be convinced that considering the initial error lying in a compact set centered at $0$ is a suitable choice.

\subsection{Non contractive case}
If the hypotheses of Theorem \ref{ERGT} are not fulfilled, we can prove the following
\begin{theorem}\label{NonContrTheor}
For any initial state $D_0=\alpha$, there exists a unique probability measure $\phi(\alpha,\cdot)$ such that
\begin{equation}
    \label{NonContr}
     \lim_{K\to\infty}\d(\RU,\RUE)=\phi g(\alpha)
  \end{equation}
where $\phi g(\alpha)=\int g(z)\phi(\alpha,\mathrm{d} z) $.
\end{theorem}
We recall that although this result holds also for the contractive case, the IRF argument is
preferable in that case since the convergence value is independent from $\alpha$.
\subsubsection{Proof of Theorem \ref{NonContrTheor}}
Let $A_{x,n}$ the set of the points that $D_k$ can reach in $n$ steps starting from $x$, i.e., $A_{x,n}=\{\q^n x +\w \sum_{i=0}^{n-1}\alpha_i\q^i,~\alpha_i\in\{-1,0,1\} \}$
\begin{lemma}\label{equilip}
For any $f\in l_f \text{-Lip}(\ms) $, there exists a positive constant $\mathrm{M}_f$ such that $P^n f \in  \frac{\mathrm{M}_f}{1-\q}\text{-Lip}(\ms)$ for any $n\in\N$. 
\end{lemma}
\proof
We have
\begin{equation}
\max_{x\in\ms}\max_{y\in A_{x,1}} \left|\frac{\mathrm{d}}{\mathrm{d}x}P(x,y)\right| \leq \frac{c\q}{\sigma\sqrt{2\pi}}
\end{equation}
then the three functions $P(x,\q x)$, $P(x, \q x+\w)$ and $P(x,\q x-\w)$ are Lipschitz with constant $l:=\frac{c\q}{\sigma\sqrt{2\pi}}$
in $\ms$. For any $x,x_0 \in \ms$ and any $(y,y_0)\in\{(\q x,\q x_0)$, $(\q x+\w, \q x_0+\w), (\q x-\w, \q x_0-\w) \}$,
\begin{equation}
  \begin{split}
&\left|P(x,y)f(y)-P(x_0,y_0)f(y_0) \right|=\\ 
&~~~~=\left|P(x,y)f(y)\pm  P(x_0,y_0)f(y) - P(x_0,y_0)f(y_0) \right| \\  
&~~~~\leq ||f||_{\infty}|P(x,y)-P(x_0,y_0)|+ P(x_0,y_0)|f(y) -f(y_0)| \\  
&~~~~\leq ||f||_{\infty} l |x-x_0|+ P(x_0,y_0)l_f|y-y_0| \\  
&~~~~\leq \left(||f||_{\infty} l+ P(x_0,y_0)l_f \q\right) |x-x_0|.  
  \end{split}
\end{equation}
Thus, 
\begin{equation}
  \begin{split}
&\left|Pf(x)-Pf(x_0) \right|=\\&=|P(x,\q x)f(\q x)-P(x_0,\q x_0)f(\q x_0)+ P(x,\q x+\w)f(\q x+\w)\\&~~~~-P(x_0,\q x_0+\w)f(\q
x_0+\w)+P(x,\q x-\w)f(\q x-\w)\\&~~~~-P(x_0,\q x_0-\w)f(\q x_0-\w)|  \\
&\leq \left(3 ||f||_{\infty} l |x-x_0|+l_f \q\right) |x-x_0|.    
  \end{split}
\end{equation}
In conclusion, $Pf\in L_1\text{-Lip}(\ms)$  where $L_1=3 l||f||_{\infty} +l_f \q$. 

Now, given any $n\in\N$ and $x=x_0+\delta$,
 \begin{equation}\label{Pngeneral}
  \begin{split}
&\left|P^nf(x)-P^n f(x_0) \right|= \left|\sum_{z\in A_{x,1}}P(x,z)  P^{n-1}f(z)- \sum_{z_0\in A_{x_0,1}}P(x_0,z_0)  P^{n-1}f(z_0)
\right| \\
&=\bigg|\sum_{z\in A_{x,1}}P(x,z)  P^{n-1}f(z) \pm \sum_{z_0\in A_{x_0,1}}P(x_0,z_0)P^{n-1}f(z_0+\q \delta)\\&~~~~ - \sum_{z_0\in
A_{x_0,1}}P(x_0,z_0)  P^{n-1}f(z_0) \bigg| \\
&= \bigg|\sum_{z\in A_{x,1}} \left[P(x,z) - P(x-\delta,z-\q \delta)\right]  P^{n-1}f(z)\\&~~~~+ \sum_{z_0\in A_{x_0,1}}P(x_0,z_0)\left[
P^{n-1}f(z_0+\q \delta)  -  P^{n-1}f(z_0)\right] \bigg| \\
&\leq \sum_{z\in A_{x,1}} \left|P(x,z) - P(x-\delta,z-\q \delta)\right|~||f||_{\infty}\\&~~~~+ \sum_{z_0\in A_{x_0,1}}P(x_0,z_0)\left|
P^{n-1}f(z_0+\q \delta)  -  P^{n-1}f(z_0)\right| \\
&\leq 3 l\delta ||f||_{\infty}+  \sum_{z_0\in A_{x_0,1}}P(x_0,z_0)\left| P^{n-1}f(z_0+\q \delta)  -  P^{n-1}f(z_0)\right|.
 \end{split}
\end{equation}
If $n=2$,
\begin{equation*}
\left|P^2f(x)-P^2 f(x_0) \right|\leq 3 l \delta ||f||_{\infty}+  \sum_{z_0\in A_{x_0,1}}P(x_0,z_0)L_1\q\delta\leq  \left(3 l ||f||_{\infty}+L_1\q\right)\delta 
\end{equation*}
that is,  $P^2f\in L_2\text{-Lip}(\ms)$ where $L_2=3 l ||f||_{\infty}+L_1\q$.
At this point, by iterating \eqref{Pngeneral}, we obtain that for any $n\in\N$, $P^nf\in L_n\text{-Lip}(\ms)$ where $L_n=3 l ||f||_{\infty}+L_{n-1}\q$.
Moreover, be recursion, 
\begin{equation*}
  L_n=3 l ||f||_{\infty}(1+\q+\dots \q^{n-1})+\q^n l_f \leq \frac{\mathrm{M}_f}{1-\q}~~~~~\mathrm{M}_f:=\max\{3 l  ||f||_{\infty}, l_f \}.
\end{equation*}
\qed

Let us recall that a sequence of measures $\{\mu_n\}_{n\in\N}$ is said to be \emph{weakly convergent} to a measure $\mu$ if $\lim_{n\to\infty}\int f(x)\mathrm{d} \mu_n=\int f(x)\mathrm{d}\mu$ for every continuous and bounded function $f$. In the next, we will denote weakly convergence by $\mu_n\xrightarrow{~w~}\mu$.

\begin{lemma}\label{exist_subs}
Let $\bP_N(x,\cdot)=\frac{1}{N}\sum_{n=0}^{N-1}P^n(x,\cdot)$, $N\in\N$. For any $x\in\ms$, there exist  a subsequence $\bP_{N_j}(x,\cdot)$, $j, N_j\in\N$, and a probability measure $\phi(x,\cdot)$ such that $\bP_{N_j}(x,\cdot)\xrightarrow{~w~} \phi(x,\cdot)$. 
\end{lemma}
\proof
This is a simple consequence of Prohorov's Theorem (see, e.g., \cite[Theorem 6.1]{bill68}; in our context tightness is trivial since the space $\ms$ is compact).
\qed

\begin{lemma}\label{all_subs}
 If all the convergent subsequences of $\bP_{N}(x,\cdot)$ weakly converge to the same $\phi(x,\cdot)$, then also $\bP_{N}(x,\cdot)$ weakly converges to $\phi(x,\cdot)$.
\end{lemma}
\proof
Again this is a consequence of Prohorov's Theorem (see, e.g.,\cite[Theorem 2.3]{bill68} )
\qed

Given  Lemmas \ref{exist_subs} and \ref{all_subs}, to prove Theorem \ref{NonContrTheor} it is sufficient to show that all the convergent subsequences of $\bP_{N}(x,\cdot)$ converge to $\phi(x,\cdot)$. 

Let us suppose that there exist a subsequence $\{M_i\}_{i\in \N} \neq \{N_j\}_{j\in\N}$ and a probability measure $\psi(x,\cdot)\neq \phi(x,\cdot)$ on $\ms$ such that  $\bP_{M_i}(x,\cdot )\xrightarrow{~w~} \psi(x,\cdot)$.

First notice that for any $m\in\N$, by applying the  Dominated Convergence Theorem,
\begin{equation*}
  \begin{split}
P^m \phi f(x)&=\int_{y\in\ms}P^m(x,\mathrm{d}y)\phi f(y)\\&=\int_{y\in\ms}P^m(x,\mathrm{d}y)\lim_{j\to +\infty}\int_{z\in\ms}\frac{1}{N_j}\sum_{n=0}^{N_j-1}P^n(y,\mathrm{d}z)f(z)\\
&=\lim_{j\to +\infty}\frac{1}{N_j}\sum_{n=0}^{N_j-1}P^{n+m}f(x)\\&=\lim_{j\to +\infty}\frac{1}{N_j}\left(\sum_{n=0}^{N_j-1}P^{n}f(x)+\sum_{n=N_j}^{N_j-1+m}P^{n}f(x)-\sum_{n=0}^{m-1}P^{n}f(x)\right)\\&=\phi f(x).\\
\end{split}
\end{equation*}

Similarly, exploiting the continuity of  $P^mf$, we obtain
$$\phi P^m f(x)=\phi f (x).$$
The same can be clearly said for $\psi$.

Now, for any $M_i\in\N$
$$\bP_{M_i}\phi f(x)=\phi f(x)$$ and $$\lim_{i\to\infty}\bP_{M_i}\phi f(x)=\phi f(x).$$
If $f$ is $l_f$-Lip$(\ms)$, then $\bP_{N_j}f(x)$ are equicontinuous by Lemma \ref{equilip}, and clearly also equibounded by $||f||_{\infty}$. Therefore, by Ascoli-Arzelà Theorem, $\phi f(x)$ is continuous and $\lim_{i\to\infty}\bP_{M_i}\phi f(x)=\psi \phi f(x)$
In conclusion,

$$\phi f (x)=\psi \phi f (x).$$
Analogously, one can prove that $\psi f(x)=\phi \psi f(x)$ and, since by Dominated Convergence $\psi \phi f(x)=\phi \psi f(x), $ we  obtain $\phi f(x)=\psi f(x)$. To summarize, we have proved that, for any $x\in \ms$, there exists a unique probability measure $\phi(x,\cdot)$ such that
$$\lim_{N\to\infty}\bP_N f(x)=\phi f(x)$$ for any $f\in l_f$-Lip($\ms$).

The thesis of Theorem \ref{NonContrTheor} follows by considering $f=g$.\qed 

The arguments used in this proof partially arise from the proof of \cite[Theorem 12.4.1]{mey93}.
  
\subsection{Simulations vs Theoretical Results}
The convergence values $\int_{\mathcal{D}}g \mathrm  d\mu$ and $\phi g(\alpha)$ can be numerically evaluated. \begin{figure}
   \begin{center}
   \includegraphics[width=5.7cm]{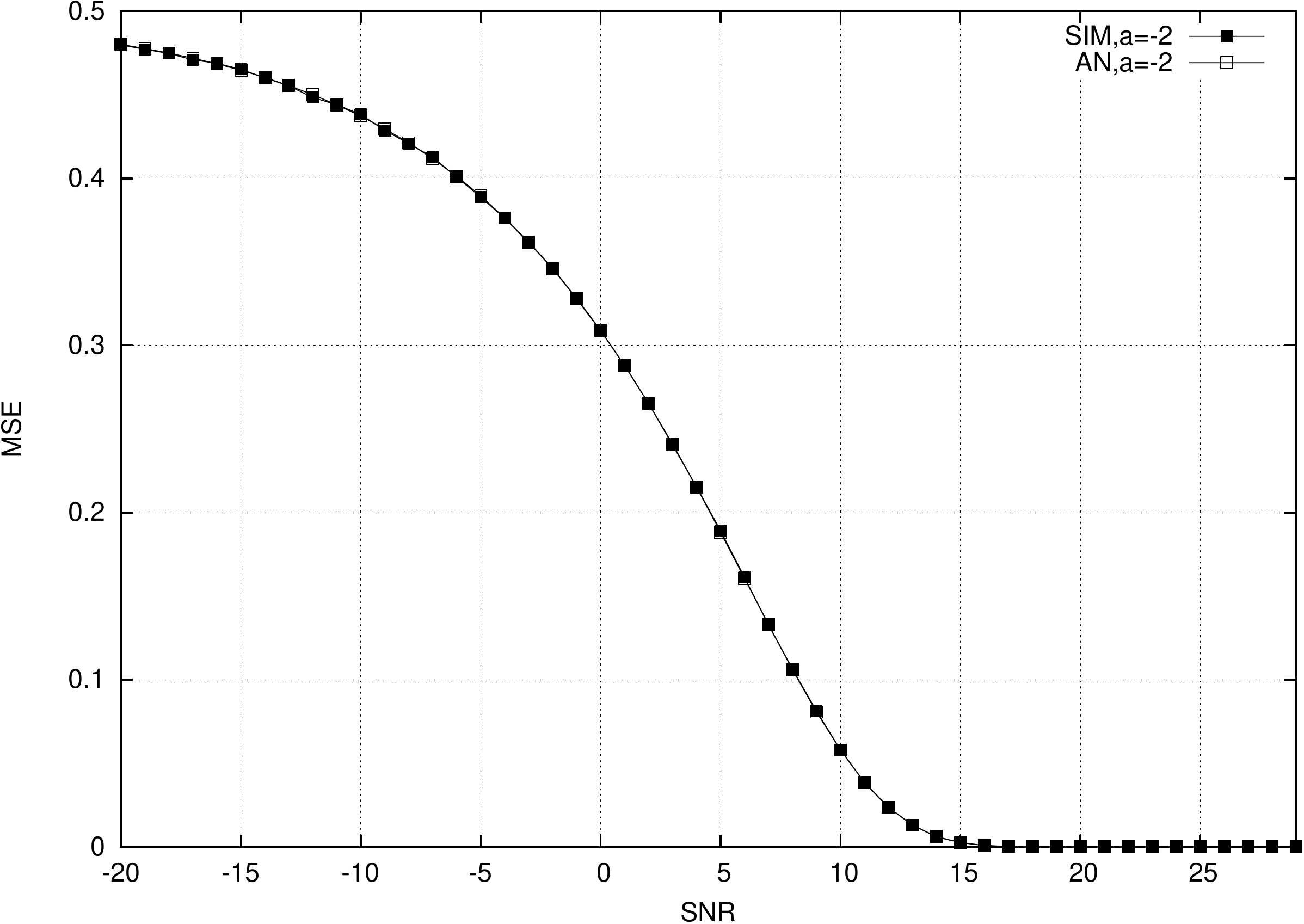}
\includegraphics[width=5.7cm]{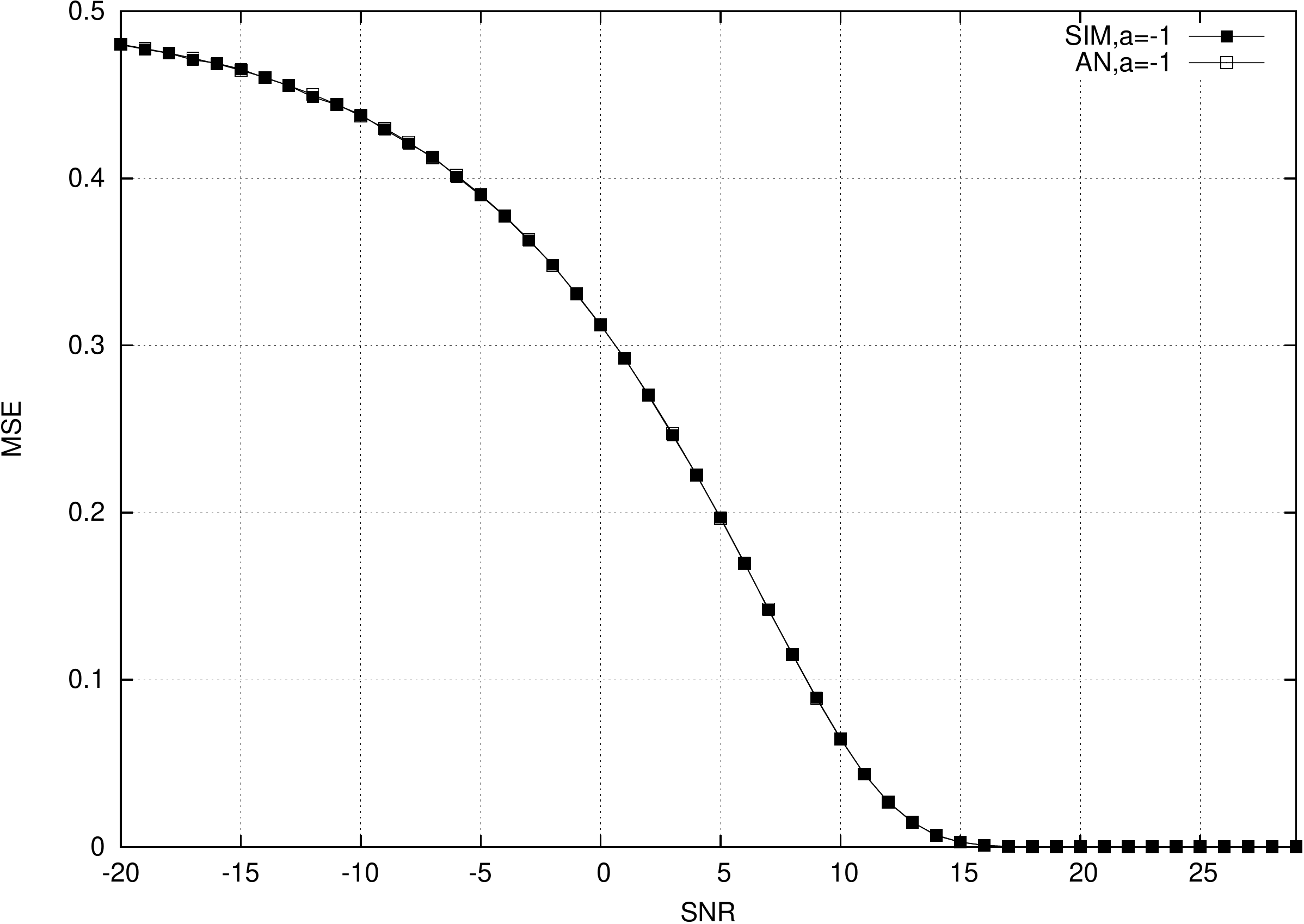}
\includegraphics[width=5.7cm]{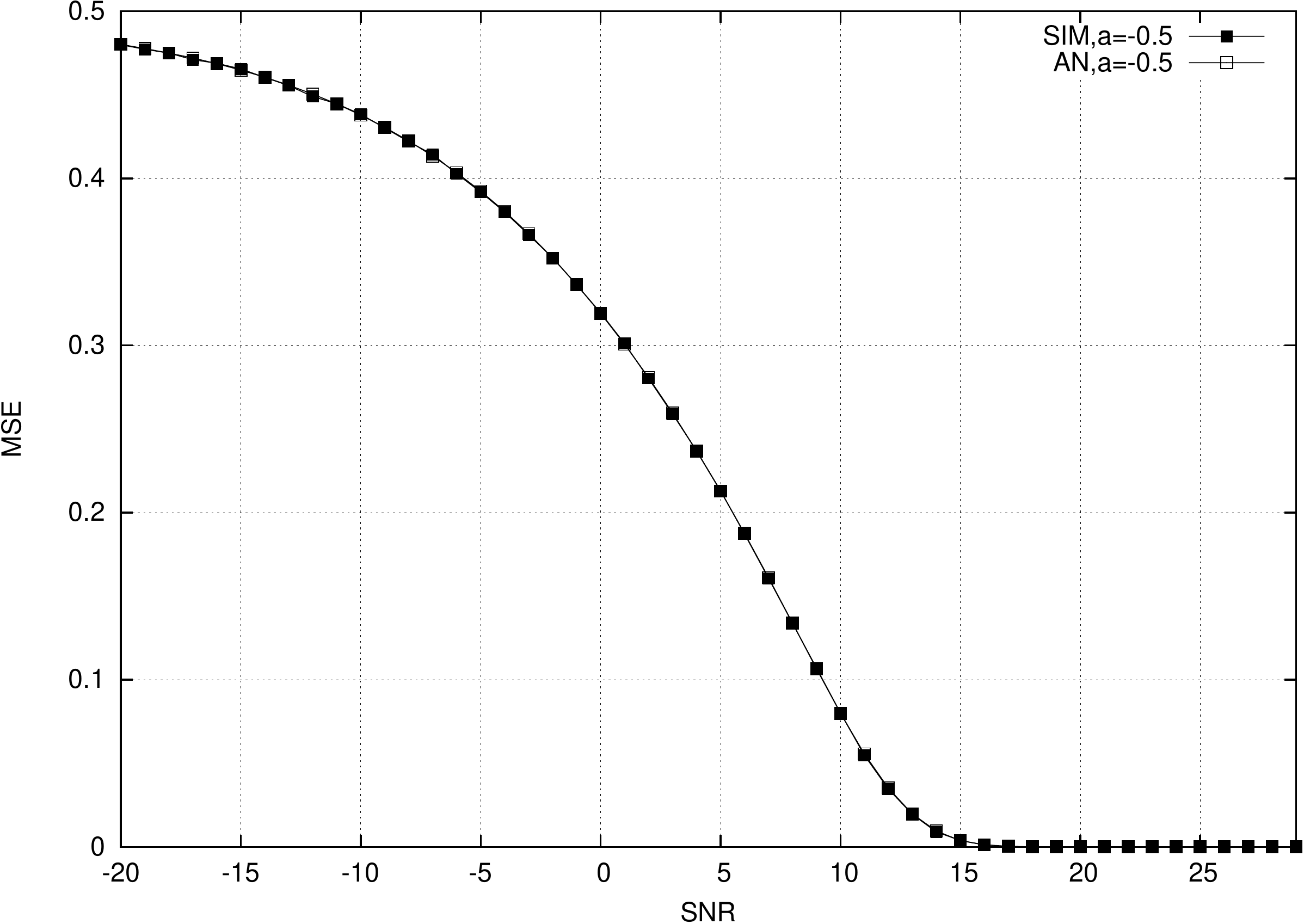}
\includegraphics[width=5.7cm]{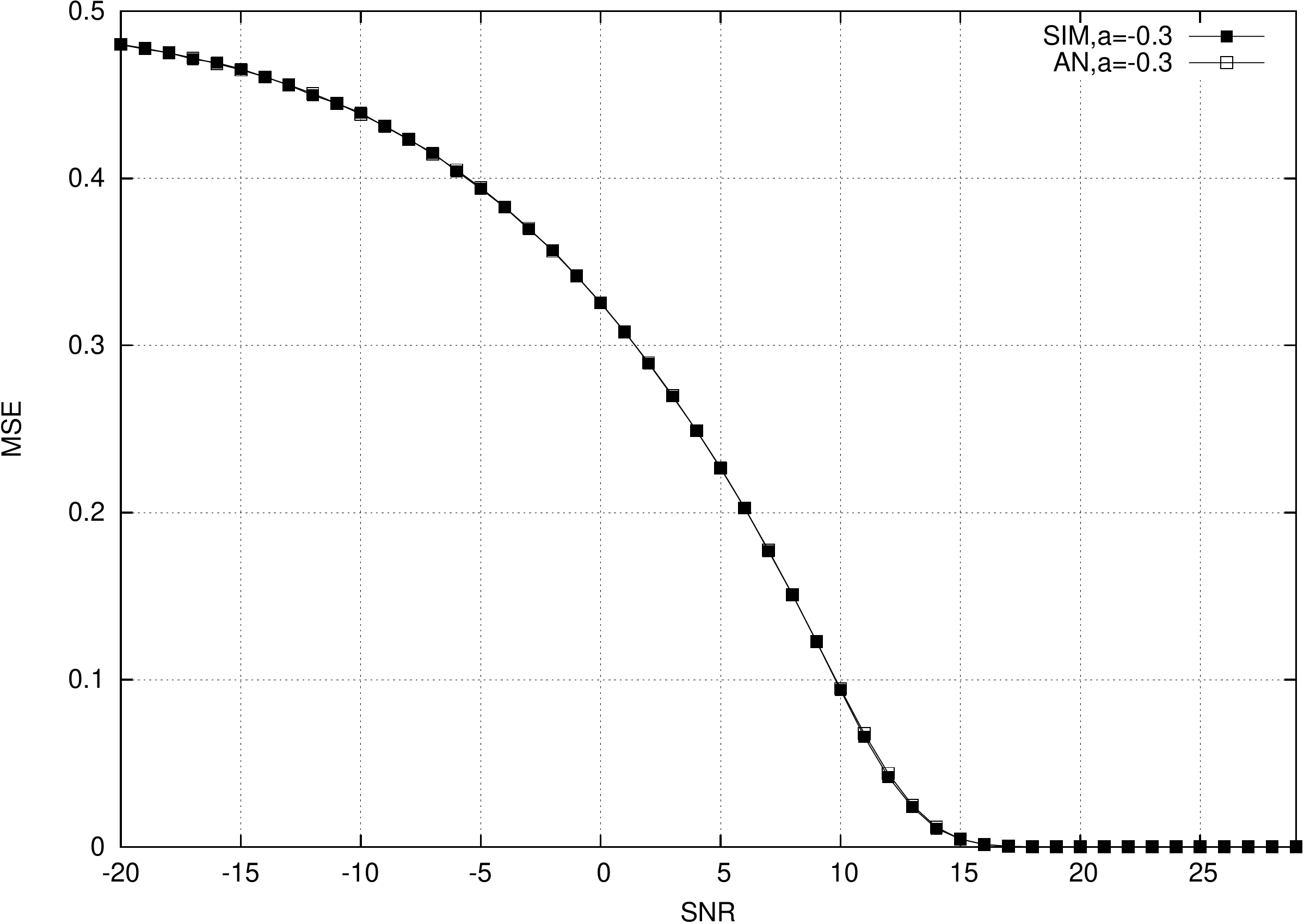}
   \caption{Simulations vs Theoretical Results: MSE for $b=c=1$,
$a=-2,-1,-0.5,-0.3$.}\label{F2}
   \end{center}
 \end{figure} In Figure \ref{F2}, we show  their consistency with the simulations previously presented. Notice that for simulations we have assumed to know the initial state $x_0$, so that $D_0=0$. Since analytical results are asymptotic,
while simulations' results are obtained by averaging  transmissions of 320 bits, we intuitively conclude that the rate of
convergence is fast. 
\section{Conclusion}\label{s_concl}
In this paper, we have proposed using the One State decoding algorithm to recover the binary input of a linear system and  we have
analyzed its behavior. When the system has particular contractive properties, the analysis is based on Iterated Random Functions, while
in the non-contractive case known results from convergence of probability measures can be exploited. The theoretical results allow to predict the Mean Square
Error of the One State Algorithm for long-time transmissions, given the parameters of the system and some prior probabilistic information. Simulations and theoretical results are consistent.

The One State Algorithm could be extended to multi-dimensional problems and to the recovery of digital inputs arising from larger source alphabets and with different probabilistic distributions. Moreover, its use for problems with feedback, such as channel equalization, should be further studied.  

\section{Acknowledgements}
The author wishes to thank Prof. Fabio Fagnani who suggested the problem and the possible solutions.





\bibliographystyle{model1-num-names}
\bibliography{references}







\end{document}